\documentclass[10pt]{article}
\usepackage{amsmath,amsthm,amsfonts}
\usepackage{graphicx}
\usepackage{hyperref}
\usepackage[usenames]{color}

\newtheorem{theorem}{Theorem}[section]
\newtheorem{lemma}[theorem]{Lemma}
\theoremstyle{definition}

\newtheorem{corollary}[theorem]{Corollary}

\newcounter{comcount}

\title{On the Reidemeister spectrum and the $R_{\infty }$ property for some free nilpotent groups}

\author{E.G. KUKINA AND V. ROMAN'KOV\footnote{The  authors were partially supported by RFBR, Grant 07-01-00392}
}

\begin{document}

\maketitle

\begin{abstract}
We describe the Reidemeister spectrum $Spec_RG$ for $G = N_{rc},$
the free nilpotent group of rank $r$ and class $c,$ in the cases:
$r \in {\mathbf N}$ and $c = 1;$ $r = 2, 3$ and $c = 2;$ $r =2$
and $c = 3,$ and prove that any group $N_{2c}$ for $c \geq 4$
satisfies to the $R_{\infty }$ property. As a consequence we
obtain that every free solvable group $S_{2t}$ of rank $2$ and
class $t \geq 2$ (in particular the free metabelian group $M_2 =
S_{22}$ of rank $2$) satisfies to the $R_{\infty }$ property.
Moreover, we prove that any free solvable group $S_{rt}$ of rank
$r \geq 2$ and class $t$ big enough also satisfies to the
$R_{\infty }$ property.{}{}
\end{abstract}

\tableofcontents
\bigskip

\section{Introduction}
\label{se:intro}

\bigskip Let $G$ be a group, and $\varphi : G \rightarrow G$ be an
automorphism of $G.$ One says that the elements $g, f \in G$ are
$\varphi -${\it twisted conjugated}, denoted by $g \sim_{\varphi }
f,$ if and only if there exists $x \in G$ such that  $(x\varphi )g
= fx.$ A class of equivalence $[g]_{\varphi }$ is called the {\it
Reidemeister class} (or the $\varphi -${\it conjugacy class of }
$\varphi $). The number $R(\varphi )$ of Reidemeister classes is
called the {\it Reidemeister number of} $\varphi .$ We define the
{\it Reidemeister spectrum} of $G$ by
\begin{equation}
\label{eq:1} Spec_R(G) = \{ R(\varphi )| \varphi \in Aut(G)\}.
\end{equation}

One says that a group $G$ {\it has the $R_{\infty }$ property for
automorphisms}, denoted by $G \in R_{\infty },$ if for every
automorphism $\varphi : G \rightarrow G$ one has $R(\varphi ) =
\infty .$ The class of groups with $R_{\infty }$ property is very
interesting in view of various of applications in
Nielsen-Reidemeister fixed point theory, representation theory,
dynamic systems, algebraic geometry and so on. See for instances
\cite{FG1}, \cite{FG2}, \cite{FGW}, \cite{FLT}, \cite{GW1},
\cite{GW2}, \cite{AG}, \cite{Sel}, \cite{TW}, \cite{TWo},
\cite{TWo2}.

In the paper by D. Gon\c{c}alves and P. Wong \cite{GW2} mainly
devoted to finitely generated nilpotent groups it was shown that
any free nilpotent group $N_{2c}$ of rank $2$ and class $c \geq 8$
satisfies to $R_{\infty }$ property. On the other hand the authors
of \cite{GW2} noted that they do not know how to extend their
techniques to $N_{rc}$ for $r \geq 3,$ and $c \geq 3.$ In the
paper by V. Roman'kov \cite{R} it was proved that any free
nilpotent group $N_{rc}$ of rank $r = 2$ or $r = 3$ and class $c
\geq 4r,$ or rank $r \geq 4$ and class $c \geq 2r,$ satisfies to
the $R_{\infty }$ property. This result provides a purely
algebraic proof of the fact that any absolutely free group $F_r, \
r \geq 2,$ has the $R_{\infty }$ property. Note that in \cite{LL}
this statement was derived by group-geometrical techniques.

In this paper we mainly consider the free nilpotent groups of
small class. We obtain the Reidemeister spectrum of a group
$N_{rc}$ for $r \in {\mathbf N}$ and $c = 1$ (free abelian case),
for $r = 2, 3$ and $c = 2$ ($2-$nilpotent case),  for $r = 2$ and
$c = 3$ ($3-$nilpotent case). As a main statement we prove that
any group $N_{2c}$ for $c \geq 4$ satisfies to the $R_{\infty }$
property. This result completes the case of rank $2.$ As a
consequence we obtain that every free solvable group $S_{2t}$ of
rank $2$ and class $t \geq 2$ (in particular the free metabelian
group $M_2 = S_{22}$ of rank $2$) satisfies to the $R_{\infty }$
property. Moreover every free solvable group $S_{rt}$ of rank $r
\geq 3$ and class $t$ big enough also satisfies to the $R_{\infty
}$ property.

\bigskip

\section{Preliminaries}
\label{se:prelim}

\bigskip

The results of this section in fact are known by paper \cite{GW2}.
Nevertheless we present them for completeness of the paper.

 Let $G$ be a finitely generated group, and let $C$
be a central subgroup of $G.$ For any automorphism $\varphi : G
\rightarrow G$ define a central subgroup
\begin{equation}
\label{eq:2} L(C, \varphi ) = \{ c \in C | \exists x \in G :
x\varphi = cx\}.
\end{equation}

It was shown in \cite{R} that any pair of elements $c_1, c_2 \in
C$ are $\varphi -$conjugated in $G$ if and only if

\begin{equation}
\label{eq:3} c_1^{-1}c_2 \in L(C, \varphi ).
\end{equation}

Thus in the case of a finitely generated abelian group $A$ the set
of all $\varphi -$conjugacy classes coincides with the set of all
cosets $A$ w.r.t. $L(A, \varphi ).$ Note that
\begin{equation}
\label{eq:4} L(A, \varphi ) = Im (\varphi - id ).
\end{equation}

Hence
\begin{equation}
\label{eq:5} R(\varphi ) = [A : L(A, \varphi )].
\end{equation}

Let $A(r) = {\mathbb Z}^r$ be a free abelian group of rank $r \in
{\mathbb N},$ and $\varphi : A(r) \rightarrow A(r)$ be any
automorphism.

Easy to see that
\begin{equation}
\label{eq:6} Spec({\mathbb Z}) = \{2\} \cup \{ \infty \}.
\end{equation}

We claim that for $r \geq 2$ the spectrum is full, i.e.

\begin{equation}
\label{eq:7}
 SpecA(r) = {\mathbb N} \cup \{\infty \}.
\end{equation}

To prove it enough to find for every number $k \in {\mathbb N}$ a
matrix $A(r, k) \in GL_r({\mathbb Z})$ such that $rank(A(r, k) -
E) = k.$

If $r = 2$ or $r = 3$ one can take
\begin{equation}
\label{eq:8} A(2, k) = \left( \begin{array}{cc} - k & 1\\
1 & 0\\
\end{array} \right),
\end{equation}

\noindent and

\begin{equation}
\label{eq:9} A(3, k) = \left( \begin{array}{ccc}
1 &  k & 1\\
1 & 1& 0 \\
1 & 0 & 0\\
\end{array} \right),
\end{equation}

\noindent respectively.

In general case for odd $r = 2t + 1$ we put

\begin{equation}
\label{eq:10} A(r, k) = diag(A(3, k), A(2, 1), ..., A(2, 1)),
\end{equation}

\noindent and for even $r = 2t$

\begin{equation}
\label{eq:11} A(r, k) = diag( A(2, k), A(2, 1), ..., A(2, 1)),
\end{equation}

\noindent where $A(2, 1)$ repeats $t - 1$ times.

It remains to note that $R(id ) = \infty $ for any group $A(r).$

Let $N$ be a finitely generated torsion free nilpotent group of
class $k.$ Let

\begin{equation}
\label{eq:12} \zeta_0N = 1 < \zeta_1N < ... < \zeta_{k-1}N <
\zeta_kN = N
\end{equation}

\noindent be the upper central series in $N.$ It is well known
(see \cite{H}, \cite{KRRRC}) that all quotients

\begin{equation}
\label{eq:13} N_i = N/\zeta_iN, \  A_i = \zeta_{i+1}N/\zeta_iN, \
i = 0,1, ..., k-1, \end{equation}

\noindent
 are finitely generated torsion free groups. In
particular, every $A_i,  \  i = 0, 1, ..., k-1,$ is a free abelian
group of a finite rank.

Let $\varphi : N \rightarrow N$ be any automorphism. Then there
are the induced automorphisms

\begin{equation}
\label{eq:14} \varphi_i : N_i \rightarrow N_i, \  \bar{\varphi_i}
: A_i \rightarrow A_i, \  i = 0, 1, ..., k-1.
\end{equation}

If $R(\varphi_i) =  \infty $ for some $ i = 0, 1, ..., k-1,$ then
$R(\varphi_j) = \infty $ for every $j < i,$ in particular,
$R(\varphi ) = \infty $ in $N.$ Moreover, if for some $ i = 0, 1,
..., k-1,$ there is a non trivial  element $\bar{a} \in A_i$ such
that $\bar{a}\bar{\varphi_i} = \bar{a},$ then $R(\varphi_j) =
\infty $ for every $j < i,$ and again $R(\varphi ) = \infty .$ To
explain the last assertion we assume that $i$ is maximal with the
property $\bar{a}\bar{\varphi_i} = \bar{a} \not= 1.$ Then the
group $N_{i+1}$ does not admit a non trivial element $x$ such that
$x\varphi_{i+1} = x.$ Hence we have $L(N_i, \varphi_i) = L(A_i,
\bar{\varphi_i}).$ Obviously $[A_i : L(A_i, \bar{\varphi_i})] =
\infty .$ Then by Lemma 2.1 from \cite{R} we derive that
$R(\varphi_i) = \infty $ in $N_i,$ and so $R(\varphi_j) = \infty $
for any $j < i,$ in particular $R(\varphi ) = \infty .$

Suppose that $R(\varphi ) < \infty .$ Then we have that the
subgroup of $\varphi_i -$invariant elements $Fix_{\varphi_i}(N_i)
= 1$ for all $i = 0, 1, ..., k-1.$ We have $[A_i : L(N_i,
\varphi_i)] = [A_i : L(A_i, \bar{\varphi_i})] = q_i,  \  i = 0, 1,
..., k-1.$ Then we have

\begin{lemma}
\label{le:2.1} Let $N$ be a finitely generated torsion free
nilpotent group of class $k,$ and $\varphi : N \rightarrow N$ be
any automorphism.  Suppose that  $R(\varphi ) < \infty ,$ and in
notions as above $[A_i : L(N_i, \varphi_i)] = [A_i : L(A_i,
\bar{\varphi_i})] = q_i, \  i = 0, 1, ..., k-1.$ Then
\begin{equation}
\label{eq:15} R(\varphi ) = \prod_{i=0}^{k-1}q_i.
\end{equation}
\end{lemma}

 Proof. The formula (\ref{eq:15}) is based on Lemma 2.4 from
 \cite{R}. In this lemma a group $G$ is considered with a central
 $\varphi -$admissible subgroup $C$ for an automorphism $\varphi :
 G \rightarrow G.$ It states that the pre image of any
 $\bar{\varphi }-$conjugacy class $[\bar{g}]_{\bar{\varphi }}$ of
 the induced automorphism $\bar{\varphi } : G/C \rightarrow G/C$
 is a disjoint union of $s = [C : L(C, \varphi_g)]$ $\varphi
 -$conjugacy classes. Here $\varphi_g = \varphi \circ \sigma_g,$
 where $\sigma_g \in InnG, \sigma_g : h \mapsto g^{-1}hg$ for all
 $h \in G.$

 We apply this statement consequently to groups $G = N_i$ and
 central subgroups $C = A_i, i = 0, 1, ..., k-1.$ By our
 assumption $Fix_{\varphi_{i+1}}(N_{i+1}) = 1,$ so if $x\varphi_i
 = cx, c \in A_i,$ then $x \in A_i.$ We see that $L(A_i, \varphi_i
 ) = L(A_i, (\varphi_i)_g)$ for every $g \in N_i.$ Hence every pre
 image of $[\bar{g}]_{\varphi_{i+1}}$ in $N_i$ is a disjoint
 union of exactly (independent of $g$) $q_i = [A_i : L(A_i, \bar{\varphi_i})]$
 $\varphi_i-$conjugacy classes. It follows that

 \begin{equation}
 \label{eq:16}
 R(\varphi_i) = R(\varphi_{i+1})\cdot q_i, \  i = 0, 1, ..., k-1.
 \end{equation}

 Repeating such process we derive (\ref{eq:15}).

\bigskip

\section{The Reidemeister specrum of $N_{rc}$ for $r = 2$ and $c = 2, 3;$ and $r = 3$ and $c = 2$}
\label{se:respec}

\bigskip
Let $N = N_{22}$ be the free nilpotent group of rank $2$ and class
$2$ (also known as the discrete Heisenberg group). Let $x, y $ be
a free basis of $N.$ Then the center $\zeta_1N$ (which coincides
with the derived subgroup $N'$) is generated by a single basic
commutator $(x, y).$

Let an automorphism $\varphi : N \rightarrow N$ induces the
automorphism of abelianization $\bar{\varphi } : N/N' \rightarrow
N/N',$ with matrix

\begin{equation}
 \label{eq:17}
 A = \left( \begin{array}{cc}
 \alpha & \beta\\
 \gamma & \delta \\
 \end{array}
 \right),
 \end{equation}

 \noindent
 matching to the basis $x, y.$ We assume that $det(A - E) = trA = k
 \not= 0,$ because in other case $R(\bar{\varphi }) = R(\varphi )
 = \infty .$ Moreover, since $(x, y)\varphi = (x, y)^{det(A)}$ we assume that $det(A) = -1,$ because
 $det(A) = 1$ implies $(x, y)\varphi = (x, y)$ and so $R(\varphi )
 = \infty $ again. Now we have $(x, y)\varphi = (x, y)^{-1},$ and so
 $[\zeta_1N : L(N, \varphi )] = [ \zeta_1N : L(\zeta_1N, \varphi )]
 = 2. $ By Lemma \ref{le:2.1}  we obtain that $R(\varphi )
 = 2|trA| \in 2{\mathbb N}.$ It remains to take the matrices $A(2,
 k)$ from the previous section to conclude that

 \begin{equation}
 \label{eq:18}
 Spec_R(N_{22}) = 2{\mathbb N} \cup \{\infty \}.
 \end{equation}

In paper \cite{Ind} it was proved that each even number belongs to
$Spec_R(N_{22}),$ but nothing was said about odd numbers.

Let now $N = N_{23}$ be the free nilpotent group of rank $2$ and
class $3.$ Let $x, y$ be a basis of $N.$

Then the center $\zeta_1N$ has a basis consisting from the basic
commutators of weight $3:$ $g_1 = (x, y, x), g_2 = (x, y, y).$
Here and so far we suppose that the brackets in any long
commutator stand from left to right, in particular $(h_1, h_2 ,
h_3) = ((h_1, h_2), h_3),$ and so on.

It is well known (see \cite{H}) that every automorphism of a free
abelian quotient  $N/N'$ of any free nilpotent group $N$ is
induced by some automorphism of $N$ itself. Moreover any
endomorphism $\eta : N \rightarrow N,$ invertible $mod(N')$
(inducing an automorphism in the abelianization $N/N'$), is
automorphism. We will use this fact later many times without
mention.

 By direct calculation we derive that an automorphism $\varphi :
N \rightarrow N,$ with matrix  on the abelianization $N/N'$ as in
(\ref{eq:17}) (we assume again that $det(A) = -1$ and $det(A - E)
= trA \not= 0)$ induces the automorphism $\bar{\varphi_0} :
\zeta_1N \rightarrow \zeta_1N$ with matrix in the basis $g_1, g_2$

\begin{equation}
\label{eq:19} A_3 = \left(
\begin{array}{cc} -\alpha & - \beta\\
- \gamma & -\delta \\
\end{array}\right).
\end{equation}

We see that

\begin{equation}
\label{eq:20} det(A_3)  = - 1, det(A_3 - E) = \alpha + \delta =
trA.
\end{equation}

Since every number ${\mathbf N}$ is obviously realized as $k =
trA_{[k]}$ for some matrix $A_{[k]} \in GL_2({\mathbb Z})$ we
complete our consideration as above to conclude that

\begin{equation}
\label{eq:21} Spec_R(N_{23}) = \{ 2k^2 | k \in {\mathbb N}\} \cup
\{\infty \}.
\end{equation}

At last, let $N = N_{32}$ be the free nilpotent group of rank $3$
and class $2.$ Let $x, y, z$ be any basis of $N.$ Then the center
$\zeta_1N $ (which coincides with the derived subgroup $N'$) has a
basis $h_1 = (x, y), h_2 = (x, z),$ and $h_3 = (y, z).$

Let an automorphism $\varphi : N \rightarrow N$ induces the
automorphism $\bar{\varphi } : N/N' \rightarrow N/N'$  with a
matrix

\begin{equation}
\label{eq:22} A = (a_{ij}),  \  i, j = 1, 2, 3, \\
\end{equation}

\noindent in a basis of the abelianization $N/N'$ matching to $x,
y, z.$ Then the matrix of the induced automorphism
$\bar{\varphi_0} : N' \rightarrow N'$ is

\begin{equation}
\label{eq:23}
B = \left(\begin{array}{ccc} M_{33} & M_{32} & M_{31}\\
M_{23} & M_{22} & M_{21}\\
M_{13} & M_{12} & M_{11}\\
\end{array}\right),
\end{equation}

\noindent where $M_{ij}$ means the minor  deriving by  deleting
$i-$row and $j-$column of $B.$

We can assume that $det(A - E) \not= 0,$ in other case $R(\varphi
) = \infty .$ By direct calculation we get $det(B) = 1.$ Also by
direct calculation we obtain

\begin{equation}
\label{eq:24} det(A - E) = detA - (M_{11} + M_{22} + M_{33}) +
a_{11} + a_{22} + a_{33} - 1,
\end{equation}

\noindent and

\begin{equation}
\label{eq:25} det(B - E) = det(B) - det(A)(a_{11} + a_{22} +
a_{33}) + M_{11} + M_{22} + M_{33} - 1.
\end{equation}

We see that these numbers in view $det(A) = \pm 1$ and $det(B) =
1$ have the same parity. If they both are even then their product
is divided by $4,$ if  both are odd this product is odd too. It
follows that numbers $4l + 2$ can not appear.

On the other hand a matrix

\begin{equation}
\label{eq:26} D_{[n]} = \left( \begin{array}{ccc} n & 1 & 1\\
1 & 1 & 0\\
1 & 0 & 0\\
\end{array}\right)\end{equation}

\noindent gives an example of automorphism $\varphi (n) $ which
induces the automorphism of abelianization $N/N'$ with matrix
$D_{[n]}$ gives an example of $R(\varphi_n) = 2n -1.$

A matrix
\begin{equation}
\label{eq:26} F_{[n]} = \left( \begin{array}{ccc} n + 1 & 1 & 1\\
2 & 1 & 0\\
1 & 0 & 0\\
\end{array}\right)\end{equation}

\noindent in the same way presents an example of an automorphism
$\psi (n)$ for which $R(\psi (n)) = 4n.$

It follows that

\begin{equation}
\label{eq:27} Spec_RN_{23} = \{2n - 1| n \in {\mathbf N}\} \cup \{
4n | n \in {\mathbf N}\} \cup \{ \infty \}.
\end{equation}

\bigskip
\bigskip

\section{Every free nilpotent group $N_{2c}$ for $c \geq 4$
satisfies to the $R_{\infty }$ property}
 \label{se:freenilp}

 \bigskip
 The main result of this section will follow from the next
 principal statement.

 {\bf Theorem 1.} Let $N = N_{24}$ be the free nilpotent group of
 rank $2$ and class $4.$ Then every automorphism $\varphi : N
 \rightarrow N$ induces the automorphism $\bar{\varphi_0 } :
 \zeta_1N \rightarrow \zeta_1N$ such that $det(\bar{\varphi_0 } -
 E) = 0.$

 \medskip
 Proof. Let $x, y$ be a basis of $N.$ Then the basic commutators
 of the weight $4$ $f_1 = (x, y, x, x), f_2 = (x, y, y, x),$ and
 $f_3 = (x, y, y, y)$ present a basis of $C = \zeta_1N$ (see
 \cite{H}). Since $N$ is metabelian one has identity

 \begin{equation}
 \label{eq:28}
 (g, f, h_1, h_2) = (g, f, h_2, h_1).
 \end{equation}

 Now suppose that any automorphism $\varphi : N \rightarrow N$
 induces the automorphism of the abelianization $\bar{\varphi } :
 N/N' \rightarrow N/N'$ with a matrix (\ref{eq:17})
 in the basis matching to $x, y.$ So, $det(A) = \pm 1.$ We assume
 that $det(A) = -1,$ in other case $(x, y)\varphi = (x,
 y)mod\zeta_2N,$ and so $R(\varphi ) = \infty .$ Moreover, we can
 assume that $tr(A) \not= 0,$ by a similar reason (see previous
 section). Let $\bar{\varphi_0} : \zeta_1N \rightarrow \zeta_1N$ be the
  automorphism induced by $\varphi .$

 By direct calculation we define the matrix of $\bar{\varphi_0} $ in
 the basis $f_1, f_2,$ and $f_3:$

 \begin{equation}
 \label{eq:29}
 B = \left( \begin{array}{ccc}
 -\alpha ^2 & - 2\alpha \beta & \beta^2\\
 -\alpha \gamma & -\alpha \delta - \beta \gamma & - \beta \delta
 \\
 -\gamma ^2 & -2 \gamma \delta & - \delta ^2\\
 \end{array}\right).
 \end{equation}

By direct calculation we derive that $det(B) = -1,$ and $det(B -
E) = 0.$

Since by our assumptions $Fix_{\varphi_1}(N_1) = 1,$ where as
above $N_1 = N/\zeta_1N$ is the free nilpotent group $N_{23}$ of
rank $2$ and class $3$ (see previous section), we conclude that
$L(N, \varphi ) = L(\zeta_1N, \bar{\varphi_0}  )$ has infinite
index in $\zeta_1N.$ Hence by Lemma 2.4 from \cite{R} we obtain
that $R(\varphi ) = \infty .$

Theorem is proved.

\begin{corollary}
\label{co:4.1} Every group $N_{2c}$ for $c \geq 4$ satisfies to
the $R_{\infty } $ property. \end{corollary}

\bigskip

Since $N_{24}$ is metabelian and in view of Theorem 2 in \cite{R}
we immediately obtain

\bigskip
{\bf Theorem 2.} 1) Every free solvable group $S_{2t}$ of rank $2$
and class $t \geq 2$ (in particular, the free metabelian group
$M_2 = S_{22}$ of rank $2$) satisfies to the $R_{\infty }$
property.

2) Every free solvable group $S_{3t}$ of rank $3$ and class $t
\geq 4$ satisfies to the $R_{\infty }$ property.

3) Every free solvable group $S_{rt}$ of rank $r \geq4$ and class
$t \geq log_2(2r + 1)$ satisfies to the $R_{\infty }$ property.

\medskip
Proof. We have for every $t \geq 2$

\begin{equation}
\label{eq:30}   N_{24} = S_{2t}/\gamma_5S_{2t},
\end{equation}

\noindent where $\gamma_5S_{2t}$ is an automorphic admissible
subgroup of $S_{2t}.$ Any automorphism $\tilde{\varphi }: S_{2t}
\rightarrow S_{2t}$ induces the automorphism $\varphi : N
\rightarrow N$ for which $R(\varphi ) = \infty .$ Hence
$R(\tilde{\varphi }) = \infty $ too.

2) By Theorem 2 in \cite{R} one has $N_{3,12} \in R_{\infty }.$
Since for every $t \geq 4$

\begin{equation}
\label{eq:31}   N_{3,12} = S_{3t}/\gamma_{13}S_{3t},
\end{equation}

\noindent by the similar argument we conclude that $S_{3t}$
satisfies to the $R_{\infty }$ property.

3) By Theorem 2 in \cite{R} one has $N_{r,2t+1} \in R_{\infty }.$
Since for every $t \geq log_2(2r+1)$

\begin{equation}
\label{eq:32}   N_{r,2r} = S_{rt}/\gamma_{2r+1}S_{rt},
\end{equation}

\noindent we again derive $S_{rt} \in R_{\infty }.$

Theorem is proved.

\medskip
{\bf Remark.} The similar results can be proved for any varieties
of groups (not just the varieties ${\cal A}^t$ of all solvable
groups of given class $t$) which admit a natural homomorphisms
onto free nilpotent groups of class big enough.

\end{document}